\documentclass[a4paper,10pt]{article}

\title{Improved lower bound on an Euclidean Ramsey problem}
\author{Jerome Barkley}

\usepackage[dvips]{color}
\usepackage[utf8x]{inputenc}
\usepackage{amsmath}
\usepackage{amsfonts}
\usepackage{amssymb}
\usepackage{rotating}
\usepackage{graphicx}
\usepackage{psfrag}

\begin{document}

\renewcommand{\figurename}{Fig.}
% Different font in captions
\newcommand{\captionfonts}{\small\it}

\makeatletter  % Allow the use of @ in command names
\long\def\@makecaption#1#2{%
  \vskip\abovecaptionskip
  \sbox\@tempboxa{{\captionfonts #1: #2}}%
  \ifdim \wd\@tempboxa >\hsize
    {\captionfonts #1: #2\par}
  \else
    \hbox to\hsize{\hfil\box\@tempboxa\hfil}%
  \fi
  \vskip\belowcaptionskip}
\makeatother   % Cancel the effect of \makeatletter

\maketitle

\begin{abstract}
It was previously shown that any two-colour colouring of $K(C_n)$ must contain a monochromatic planar $K_4$ subgraph for $n \ge N^*$, where $6 \le N^* \le N$ and $N$ is Graham's number. The bound was later \cite{aemp} improved to $11 \le N^* \le N$. In this article, it is improved to $13 \le N^* \le N$.
\end{abstract}

\section{Introduction}

Consider an $n$-dimensional hypercube $C_n$. Consider the complete graph $K(C_n)$ connecting the vertices of the $n$-cube $C_n$, and consider a two-colour colouring of $K(C_n)$. Let $N^*$ be the smallest integer, such that any two-colour colouring of $K(C_{N^*})$ must contain a monochromatic planar $K_4$ (4 vertices, complete) subgraph.

It was shown \cite{rtfnps} that $6 \le N^* \le N$, where $N$ is Graham's number, a very large number. It was later shown \cite{aemp} that $11 \le N^*$, by constructing a colouring of $K(C_{10})$ which doesn't contain a monochromatic planar $K_4$ subgraph. According to Exoo, at least two other people have constructed unpublished colourings of $K(C_{11})$, showing that $12 \le N^*$.

\begin{figure}\centering
\label{cubes}
\includegraphics{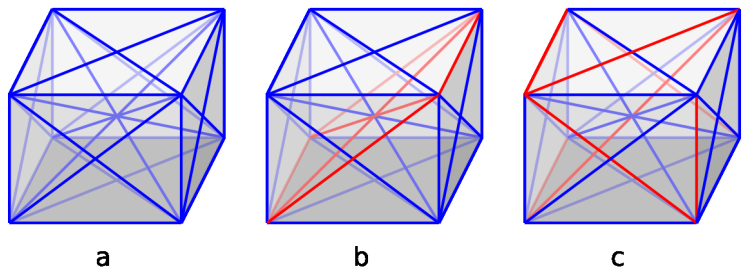}
\caption{
a: $K(C_3)$, all edges coloured blue.
b: $K(C_3)$, one planar $K_4$ subgraph coloured red.
c: $K(C_3)$, no monochromatic planar $K_4$ subgraphs.
}
\end{figure}

\textit{Fig. \!\ref{cubes}a} shows $K(C_3)$, \textit{Fig. \!\ref{cubes}b} shows $K(C_3)$ with one planar $K_4$ subgraph in red, and \textit{Fig. \!\ref{cubes}c} shows a colouring with no monochromatic planar $K_4$ subgraphs, showing that $4 \le N^*$. Colourings in higher dimensions are hard to illustrate with images.

\section{Construction of a colouring of $K(C_{10})$}

In \cite{aemp}, a colouring of $K(C_{10})$ is constructed by colouring all edges the same colour, cycling through the list of edges in some order, and flipping the colour of each edge, whenever flipping the colour of the edge would not create new monochromatic planar $K_4$ subgraphs. The algorithm terminates when there are no monochromatic planar $K_4$ subgraphs left. Such an algorithm took over 30 hours on a 1.4 GHz CPU.

In the article at hand, also, is found a colouring of $K(C_{10})$ using an almost as simple, nondeterministic, algorithm as follows. The probability function is arbitrarily chosen, perhaps almost any would do, and a better one probably exists.

Colour all edges the same colour. While there are monochromatic planar $K_4$ subgraphs remaining, randomly pick an edge $e$. Let $n_B$ be the number of monochromatic planar $K_4$ subgraphs containing $e$ and $n_G$ be the number of planar $K_4$ subgraphs containing $e$ which would become monochromatic if the edge colour was flipped. Let $P = \min(\frac{n_B}{10+100n_G}, 1)$. Flip the edge colour with probability $P$.

%r184
%As an optimisation, once only a small subset of the edges are part of monochromatic planar $K_4$ subgraphs, the implementation of this algorithm switches to randomly picking from a list of all edges which may be part of a monochromatic planar $K_4$ subgraph. This saves the implementation from constantly selecting edges which will be flipped with probability $P = 0$.

The implementation of this algorithm took 2 minutes and 39 seconds to colour $K(C_{0})$ up to $K(C_{10})$, on a single 2.4 GHz core of a Core2 Q6600 processor. It is not known how the processor or memory bandwidth compares to the processor and memory bandwidth used for the previous construction. ($K(C_{2})$ is the first non-trivial graph which actually contains (and is) a planar $K_4$ subgraph.)

\section{Construction of a colouring of $K(C_{11})$}

In \cite{aemp} was estimated that an attempt at constructing a colouring of $K(C_{11})$ would take over 1000 hours, and over half a gigabyte of memory, which was too much memory to make the attempt.

Here, the same non-deterministic algorithm was used to construct a colouring of $K(C_{11})$ as was used here for $K(C_{10})$. The implementation of this algorithm took 48 minutes and 6 seconds to colour $K(C_{11})$, on the same processor as before, using about 17MB of memory. Since the algorithm picks edges randomly, and does not need to contain an ordered list of all edges, the algorithm needs less than the previously estimated half a gigabyte of memory.

Since the algorithm terminated, $12 \le N^*$. The colouring of $K(C_{11})$ is available at http://www.nbi.dk/\~{}barkley/graham/ .
%\textcolor[rgb]{0,0,1}{\tt{\underline{http://www.nbi.dk/\~{}barkley/graham/}}}.

When running the same algorithm to colour $K(C_{12})$, it appeared that it would take a very long time.

\section{Estimates of difficulty of colouring $K(C_n)$}

The $K(C_n)$ graph has $\frac{2^n(2^n-1)}2 = 2^{2n-1}-2^{n-1} \equiv n_E$ edges.
When specifying a colouring, there are thus $n_E$ bits of freedom in specifying the colouring.

In \cite{aemp}, the $K(C_n)$ graph is shown to have $2^{n-3}(3^n-2^{n+1}+1) \equiv n_K$ planar $K_4$ subgraphs.
There are $2^6 = 64$ ways of colouring a $K_4$ graph, $62$ of which are not monochromatic. That is, $\frac{31}{32}$ of all possible colourings of a $K_4$ graph are monochromatic. The constraint that one particular $K_4$ subgraph be monochromatic thus constrains $-\log_2 \frac{31}{32} \approx 0.0458$ of the bits of freedom.
Under the very naïve assumption that each constraint, that each planar $K_4$ subgraph not be monochromatic, is independent of each other constraint, $-n_K\log_2 \frac{31}{32}$ bits of freedom would be constrained. The fraction of constrained bits is then $-\frac{n_K}{n_E}\log_2 \frac{31}{32} \equiv n_F$. (So, for example, if the assumption was correct, and $n_E$ was 10 and $n_K$ was 30\%, then 7 bits would be required to specify a particular of $2^7$ solutions.)

\begin{tabular}{|r|r|r|r|}
\hline
$n$ & $n_E$ & $n_K$ & $n_F$ \\
\hline
2 & 6 & 1 & 0.763\% \\
3 & 28 & 12 & 1.963\% \\
4 & 120 & 100 & 3.817\% \\
5 & 496 & 720 & 6.649\% \\
6 & 2016 & 4816 & 10.942\% \\
7 & 8128 & 30912 & 17.420\% \\
8 & 32640 & 193600 & 27.168\% \\
9 & 130816 & 1194240 & 41.815\% \\
10 & 523776 & 7296256 & 63.805\% \\
11 & 2096128 & 44301312 & 96.805\% \\
12 & 8386560 & 267904000 & 146.317\% \\
13 & 33550336 & 1615810560 & 220.594\% \\
14 & 134209536 & 9728413696 & 332.016\% \\
\hline
\end{tabular}

Under the assumption, the problem is thus overspecified for $n \ge 12$, making a solution for $n \ge 12$ seem unlikely.

The fraction is trivially correct for $n = 2$. Out of the $2^{28} = 268435456$ colourings for $K(C_3)$, $182596118$ of them are without monochromatic planar $K_4$ subgraphs. The correct fraction for $n = 3$ is thus actually $-\frac1{28}\log_2 \frac{182596118}{2^{28}} \approx 1.985\%$, which means there are slightly fewer solutions for $n = 3$ than if the assumption were correct. It is currently impractical to test all possible colourings of $K(C_n)$ and count the solutions, for $n \ge 4$.

\subsection{Symmetries}

A note on notation — $C_n$ is used here to refer to the $n$-element cyclic group. (As opposed to the $C_n$ in $K(C_n)$, where $C_n$ is used to refer to the $n$-dimensional cube.)

It is possible to require that the colouring be symmetric, such that the colouring does not change under some subgroup of the automorphism group of the problem. The automorphism group of the $n$-cube is the signed permutation group $C_2 \wr_n S_n$ (notation\footnote{
$A \wr_n B$ is the wreath product — the direct product $A^n$ combined with $B$ represented as a permutation of the $n$ elements of $A^n$.
}) and the colour flipping symmetry $C_2$. (This risks turning a solvable problem into an unsolvable one.)

Since using colour flipping symmetry or the signedness of the signed permutation group didn't seem to help, only subgroups of $S_n$, the (unsigned) permutation group will be considered here.

Under such a symmetry, many $K_4$ subgraphs will be congruent to each other, and the constraints imposed by them are equivalent. This means means less constraints.

If two or more edges on the same $K_4$ subgraph are congruent to each other under the symmetry, or can somehow be shown to have the same colour, the constraint is reduced from a constraint on 6 edges, to a constraint on 5 or less edges (which means harder constraints). If a constraint can be reduced in that way to a constraint on 1 edge, the edge must be both colours, which is not possible, and a solution with the given symmetry is not possible.

If there are a constraint on only edge $a$ and edge $b$ (that is, a constraint that edge $a$ and edge $b$ have opposite colour), and another constraint on only edge $b$ and edge $c$, then since both edge $a$ and edge $c$ have the opposite colour as edge $b$, edge $a$ must then have the same colour as edge $c$, so they may be considered equivalent. (And thus the two constraints become equivalent.)

If there is a constraint on a set $A$ of edges and another constraint on a set $B$ of edges, and if $A \subseteq B$, then the constraint on $A$ implies the constraint on $B$, and the constraint on $B$ is redundant.

A fraction $\frac{2^\nu-2}{2^\nu}$ of colourings (respecting the symmetry) will satisfy a constraint on $\nu$ edges. Under the naïve assumption that after applying the symmetry, the remaining constraints are (still) independent, each constraint on $\nu$ edges will constrain $-\log_2\frac{2^\nu-2}{2^\nu}$ of the bits of freedom.

\subsubsection{Some arbitrarily chosen symmetries of $K(C_{9})$}

One possible symmetry of a colouring of $K(C_{9})$ is the group $S_9$. There is a colouring of $K(C_{9})$ with $S_9$ symmetry with no monochromatic planar $K_4$ subgraphs, therefore there is a colouring (the same one) with any possible (unsigned) permutation symmetry.

\subsubsection{Some arbitrarily chosen symmetries of $K(C_{10})$}

There are no solutions for $K(C_{10})$ with $S_{10}$ (or $A_{10}$) symmetry. There is a solution with $S_5$ symmetry, where $S_5$ is represented as a primitive permutation group on 10 coordinates.

Generators for $S_5$ found with GAP\footnote{
{\tt PrimitiveGroup(10, 2);}
} are \\
$(1\ 5\ 7)(2\ 9\ 4)(3\ 8\ 10)$ and
$(1\ 8)(2\ 5\ 6\ 3)(4\ 9\ 7\ 10)$. This group has order $5! = 120$.

\subsubsection{Some arbitrarily chosen symmetries of $K(C_{11})$}

Three possible symmetries of a colouring of $K(C_{11})$ are the group $M_{11}$ represented as a permutation of 11 coordinates \cite{m11}, a Sylow 3-subgroup (GAP\footnote{
{\tt SylowSubgroup(Group((1,2,3,4,5,6,7,8,9,10,11), (1,2)), 3);} or\\
{\tt SylowSubgroup(SymmetricGroup(IsPermGroup, 11), 3);} or\\
{\tt WreathProduct(CyclicGroup(IsPermGroup, 3), CyclicGroup(IsPermGroup, 3));}
})of the permutation group $S_{11}$ (which is equivalent to $C_3 \wr_3 C_3$) and the projective special linear group $L_2(11)$ (GAP\footnote{
{\tt PrimitiveGroup(11, 5);}
}). For $M_{11}$, an almost-solution which violates just one constraint (which corresponds to multiple planar $K_4$ subgraphs) exists.

\begin{tabular}{|r|c|r|}
\hline
Group & Generators & Order \\
\hline
$M_{11}$ & $(2\ 10)(4\ 11)(5\ 7)(8\ 9)$ \quad $(1\ 4\ 3\ 8)(2\ 5\ 6\ 9)$ & $11!/7! = 7920$ \\
\hline
$Syl_3(S_{11})$ & $(3\ 9\ 7)(6\ 11\ 10)$ \quad $(3\ 6\ 4)(5\ 9\ 11)(7\ 10\ 8)$ & $3^4 = 81$ \\
\hline
$L_2(11)$ & $(1\ 5)(2\ 4)(3\ 10)(7\ 11)$ &\\& $(3\ 11\ 5)(4\ 7\ 9)(6\ 8\ 10)$ & $660$ \\
\hline
\end{tabular}

\subsubsection{Some arbitrarily chosen symmetries of $K(C_{12})$}

Five possible symmetries of a colouring of $K(C_{12})$ are the group $M_{11}$ represented as a permutation of 12 coordinates \cite{m11}, a Sylow 3-subgroup of the permutation group $S_{12}$ (GAP\footnote{
{\tt SylowSubgroup(Group((1,2,3,4,5,6,7,8,9,10,11,12), (1,2)), 3);} or\\
{\tt SylowSubgroup(SymmetricGroup(IsPermGroup, 12), 3);} or\\
{\tt DirectProduct(WreathProduct(CyclicGroup(IsPermGroup, 3), CyclicGroup(IsPermGroup, 3)), CyclicGroup(IsPermGroup, 3));}
}) = $(C_3 \wr_3 C_3) \times C_3$, $(D_4)^3$, $AGL_1(5) \times L_3(2)$ (GAP\footnote{
{\tt DirectProduct(PrimitiveGroup(5, 3), PrimitiveGroup(7, 5));}
}) and $S_3 \times S_9$. No solution exists for $M_{12}$.

\begin{tabular}{|r|c|r|}
\hline
Group & Generators & Order \\
\hline
$M_{11}$ & $(1\ 6)(2\ 9)(5\ 7)(8\ 10)$ &\\& $(1\ 6\ 7\ 4)(2\ 8)(3\ 9)(5\ 11\ 12\ 10)$ & $11!/7! = 7920$ \\
\hline
$Syl_3(S_{12})$ & $(1\ 7\ 11)(3\ 4\ 10)$ &\\& $(1\ 10\ 12)(2\ 7\ 3)(4\ 8\ 11)$ \quad $(5\ 6\ 9)$ & $3^5 = 243$ \\
\hline
$(D_4)^3$ & $(1\ 2)$ \quad $(1\ 3)(2\ 4)$ \quad $(5\ 6)$ &\\& $(5\ 7)(6\ 8)$ \quad $(9\ 10)$ \quad $(9\ 11)(10\ 12)$ & $8^3 = 512$ \\
\hline
\footnotesize $AGL_1(5) \times L_3(2)$ & \small $(2\ 3\ 4\ 5)$ \quad $(1\ 2\ 3\ 5\ 4)$ &\\& \small $(6\ 9)(11\ 12)$ \quad $(6\ 8\ 7)(9\ 12\ 10)$ & \footnotesize $20 \times 168 = 3360$ \\
\hline
\end{tabular}

\subsubsection{Some arbitrarily chosen symmetries of $K(C_{13})$}

One possible symmetry is $L_3(3)$, found with GAP\footnote{
{\tt PrimitiveGroup(13, 7);}
}. Generators are \\
$(1\ 10\ 4)(6\ 9\ 7)(8\ 12\ 13)$ and
$(1\ 3\ 2)(4\ 9\ 5)(7\ 8\ 12)(10\ 13\ 11)$. This group has order $2^4 \cdot 3^3 \cdot 13 = 5616$. A colouring exists which violates 142 constraints.

\subsubsection{Some arbitrarily chosen symmetries of $K(C_{14})$}

One possible symmetry is $S_5 \times S_9$. This group has order $5! \cdot 9! = 43545600$. A colouring exists which violates 83 constraints.

\subsubsection{Revised estimates of difficulty of colouring $K(C_{n})$}

The number of constraints with each number of edges is shown, along with the fraction of constrained bits, for the identity (same as before), and for the groups. A solution has been found for groups marked with green.

%2 & 6 & 1 & 0.763\% \\
%3 & 28 & 12 & 1.963\% \\
%4 & 120 & 100 & 3.817\% \\
%5 & 496 & 720 & 6.649\% \\
%6 & 2016 & 4816 & 10.942\% \\
%7 & 8128 & 30912 & 17.420\% \\
%8 & 32640 & 193600 & 27.168\% \\
\newcommand{\solved}[1]{\textcolor[rgb]{0,0.75,0}{#1}}
\newcommand{\impossible}[1]{\textcolor[rgb]{1,0,0}{#1}}
\begin{tabular}{|r|r|r|c|r|}
\hline
Group & $n$ & $n_E$ & \impossible{1}, 2, 3, 4, 5, 6 & $n_F$ \\
\hline
%9 & 130816 & 1194240 & 41.815\% \\
\solved{$I$} &9& 130816 & 0, 0, 0, 0, 1194240 & 41.815\% \\
%constraints = 253 (0 6 0 106 0 141), freedom = 111, duty = 32.878698, difficulty = 0.296204
\solved{$S_9$} &9& 111 & 6, 0, 106, 0, 141 & 29.620\% \\
\hline
%constraints = 324 (1 6 0 106 0 211), freedom = 142, duty = inf, difficulty = inf
$S_{10}$ &10& 142 & \impossible{1}, 6, 0, 106, 0, 211 & \impossible{$\infty$\%} \\
%constraints = 65601 (0 12 64 3090 420 62015), freedom = 5432, duty = 3513.457452, difficulty = 0.646807
\solved{$S_5$} &10& 5432 & 12, 64, 3090, 420, 62015 & 64.681\% \\
%10 & 523776 & 7296256 & 63.805\% \\
\solved{$I$} &10& 523776 & 0, 0, 0, 0, 7296256 & 63.805\% \\
\hline
%11 & 2096128 & 44301312 & 96.805\% \\
\solved{$I$} &11& 2096128 & 0, 0, 0, 0, 44301312 & 96.805\% \\
%constraints = 71435 (0 4 66 3168 1340 66857), freedom = 4034, duty = 3828.755960, difficulty = 0.949121
\solved{$L_2(11)$} &11 & 4034 & 4, 66, 3168, 1340, 66857 & 94.912\% \\
%constraints = 6348 (0 12 38 901 516 4881), freedom = 562, duty = 472.956902, difficulty = 0.841560
$M_{11}$ &11& 562 & 12, 38, 901, 516, 4881 & 84.156\% \\
%constraints = 640200 (0 0 16 168 23744 616272), freedom = 36944, duty = 30477.326076, difficulty = 0.824960
\solved{$Syl_3(S_{11})$} &11 & 36944 & 0, 16, 168, 23744, 616272 & 82.496\% \\
\hline
%constraints = 3479 (1 11 33 655 172 2607), freedom = 429, duty = inf, difficulty = inf
$M_{12}$ &12& 429 & \impossible{1}, 11, 33, 655, 172, 2607 & \impossible{$\infty$\%} \\
%12 & 8386560 & 267904000 & 146.317\% \\
\solved{$I$} &12& 8386560 & 0, 0, 0, 0, 267904000 & 146.317\% \\
%constraints = 39569 (0 17 97 3104 1801 34550), freedom = 1969, duty = 2405.436473, difficulty = 1.221654
$M_{11}$ &12& 1969 & 17, 97, 3104, 1801, 34550 & 122.165\% \\
%constraints = 341887 (0 10 100 1950 15068 324759), freedom = 14138, duty = 16705.294592, difficulty = 1.181588
\solved{$C_3 \wr_4 A_4$} &12& 14138 & \!10, 100, 1950, 15068, 324759\! & 118.159\% \\
%constraints = 1022796 (0 16 74 1985 35762 984959), freedom = 41588, duty = 48873.648092, difficulty = 1.175186
\solved{$C_3 \wr_4 C_4$} &12& 41588 & 16, 74, 1985, 35762, 984959 & 117.519\% \\
%constraints = 1362824 (0 0 24 768 49912 1312120), freedom = 55440, duty = 64905.126127, difficulty = 1.170727
\solved{$Syl_3(S_{12})$} &12& 55440 & 0, 24, 768, 49912, 1312120 & 117.073\% \\
%constraints = 1550639 (0 127 2671 98585 0 1449256), freedom = 84070, duty = 86608.752163, difficulty = 1.030198
$(D_4)^3$ &12 & 84070 & \!127, 2671, 98585, 0, 1449256\! & 103.020\% \\
%constraints = 161085 (0 52 388 12350 4858 143437), freedom = 10168, duty = 9614.470576, difficulty = 0.945562
\!\!\footnotesize{$AGL_1(5) \times L_3(2)$}\!\! &12 & 10168 & \!52, 388, 12350, 4858, 143437\! & 94.556\% \\
%constraints = 14440 (0 24 80 2894 0 11442), freedom = 2234, duty = 1138.803672, difficulty = 0.509760
$S_3 \times S_9$ &12& 2234 & 24, 80, 2894, 0, 11442 & 50.976\% \\
\hline
%13 & 33550336 & 1615810560 & 220.594\% \\
$I$ &13& \!33550336\! & 0, 0, 0, 0, 1615810560 & 220.594\% \\
%constraints = 319367 (0 34 196 11857 5340 301940), freedom = 9174, duty = 16726.710300, difficulty = 1.823273
$L_3(3)$ &13& 9174 & \!34, 196, 11857, 5340, 301940\! & 182.327\% \\
\hline
%14 & 134209536 & 9728413696 & 332.016\% \\
$I$ &14& \!\!134209536\!\! & 0, 0, 0, 0, 9728413696 & 332.016\% \\
%constraints = 65636 (0 72 292 10748 0 54524), freedom = 6256, duty = 4761.140620, difficulty = 0.761052
$S_5 \times S_9$ &14& 6256 & 72, 292, 10748, 0, 54524 & 76.105\% \\
\hline
\end{tabular}

A colouring of $K(C_{11})$ with no monochromatic planar $K_4$ subgraphs can be found in under an hour with no symmetry ($I$), or under a minute\footnote{
%prob = baseProb*badness/std::min(baseProb, std::max<uint64_t>(badness, 10 + 100*(goodness*5>badness? goodness*5-badness : 0)));
Using $P = \min\left(\frac{n_B}{10+100\max(5n_G-n_B, 0)}, 1\right)$ and writing randomly into a blacklist of 3 recently flipped edges, to avoid flipping back and forth.
} with $Syl_3(S_{11})$ symmetry, or under five seconds with $L_2(11)$ symmetry.

After applying the $S_3 \times S_9$ symmetry, the problem of finding a colouring of $K(C_{12})$ with no monochromatic planar $K_4$ subgraphs appears easier and smaller than the the problem for $K(C_{11})$ with the three mentioned symmetries. However, no solution for $K(C_{12})$ was found for this particular symmetry. The solutions found for $K(C_{12})$ in the next section have $C_3 \wr_4 C_4$ and $C_3 \wr_4 A_4$ symmetry, which according to the previous table were unlikely to have solutions.

\section{Construction of a colouring of $K(C_{12})$}

A relative probability of flipping edges, $P = \min\left(\frac{n_B}{10+100\max(5n_G-n_B, 0)}, 1\right)$, is used. This arbitrarily chosen probability seems to function better than the previous arbitrarily chosen one. Also, a blacklist of 3 recently flipped edges is used, to avoid flipping back and forth. A random entry in the blacklist is overwritten each flip.

After applying a symmetry, some edges have more constraints than other edges. Let the value $x$ of a constraint $C$ be the maximum of the numbers of constraints affecting one of edges affected by $C$. For a given cutoff $\kappa$, ignore all constraints with $x \leq \kappa$.

The following algorithm is used. Start with all edges the same colour. Let $\kappa$ be the maximum $x$. While all un-ignored constraints are satisfied, reduce $\kappa$. Flip an edge. If flipping 2\,000\,000 times without reducing $\kappa$, run 2\,000\,000 more flips with $\kappa = 0$ and repeat with $\kappa$ back to the maximum $x$.

The idea is that some constraints are ``harder'' to satisfy than other constraints, and that if trying to first satisfy the ``harder'' constraints without the ``easier'' constraints getting in the way, the ``easier'' constraints will then be easy to satisfy.

The symmetry used here is $C_3 \wr_4 C_4$. This group has order $3^4 \cdot 4 = 324$.

The implementation of this algorithm took 2 hours, 31 minutes and 39 seconds, or 81\,658\,217 edge colour flips, to colour $K(C_{12})$ with no monochromatic planar $K_4$ subgraphs, on the same processor as before. (Other attempts were running at the same time on other processors, which may have reduced memory bandwidth and increased time for this attempt.)

%r219
Since the algorithm terminated, $13 \le N^*$. The colouring of $K(C_{12})$ is available at http://www.nbi.dk/\~{}barkley/graham/ .
%\textcolor[rgb]{0,0,1}{\tt{\underline{http://www.nbi.dk/\~{}barkley/graham/}}}.

The existence of a solution with $C_3 \wr_4 C_4$ symmetry implies the existence of, and is, a solution with $I$ symmetry, since $I \leq C_3 \wr_4 C_4$.

The closest to solutions that four later (simultaneous) attempts with the same symmetry, algorithm and implementation got after about 14 hours were almost-solutions that violated 9, 1, 10 and 4 constraints. This suggests that luck was a major factor in the 2 hours, 31 minutes and 39 seconds time of the first attempt with that particular symmetry, algorithm and implementation.

With four simultaneous attempts at colouring with $C_3 \wr_4 A_4$ symmetry, one attempt found a solution after 40 hours, 20 minutes and 14 seconds, after 1\,961\,430\,488 edge colour flips, another found a solution after 37 hours, 4 minutes and 12 seconds, after 1\,969\,734\,275 edge colour flips, the third was down to two violated constraints after four days, and the last was down to one violated constraint after four days.

The existence of a solution with $C_3 \wr_4 A_4$ symmetry implies the existence of, and is, a solution with $I$ symmetry and a solution with $Syl_3(S_{12})$ symmetry, since $I \leq Syl_3(S_{12}) \leq C_3 \wr_4 A_4$.

\section{Conclusion}

It is possible to bi-colour $K(C_{12})$ with no monochromatic planar $K_4$ subgraphs, therefore $13 \le N^*$.

An argument has been given why it should not be possible to bi-colour $K(C_{n})$ for $n \gtrapprox 12$ with no monochromatic planar $K_4$ subgraphs. An argument has also been given why it should be possible to bi-colour $K(C_{12})$ with no monochromatic planar $K_4$ subgraphs, with certain symmetries. However, no colouring was found for the symmetries that seem easiest — only for a symmetry which was still estimated to be impossible. It therefore seems that no estimates were accurate.

It is not surprising that no estimates seem accurate. The estimates are based on the naïve assumtion that the constraints are independent, when the constraints are actually far from being independent. A small group of constraints together may be more or less constraining than the the small group constraints would be, had they been independent, causing an error in the estimate. Especially since the problem is highly symmetric, there seems to be no reason why the errors in an estimate would tend to cancel out — the estimates could easily be wildly off.

Maybe a better estimate is possible — the current upper bound for $N^*$ seems a bit big.

\end{document}